\documentclass[11pt]{article}
\usepackage{amsthm,amsfonts,amsmath,amssymb,amscd}
\usepackage{graphicx}
\usepackage{pdfsync}
\usepackage{appendix}

\swapnumbers
\theoremstyle{plain}

\newtheorem{thm}{THEOREM}[section]
\newtheorem{lm}[thm]{LEMMA}

\theoremstyle{definition}

\theoremstyle{remark}
                                
\newcommand{\version}{September 27, 2009}
                                
\def\d{{\rm d}}
\def\dd{{\rm d}}
\def\F{{\mathcal F}}
\def\D{{\mathcal D}}
\newcommand{\R}{{\mathord{\mathbb R}}}

\newcommand{\N}{{\mathord{\mathbb N}}}

\begin{document}

\newcommand{\lanbox}{\hfill \hbox{$\, 
\vrule height 0.25cm width 0.25cm depth 0.01cm
\,$}}

\markboth{\scriptsize{CCL \version}}{\scriptsize{CCL \version}}

\bibliographystyle{plain}
\title{{\bf HARDY-LITTLEWOOD-SOBOLEV INEQUALITIES VIA  FAST DIFFUSION FLOWS}}

\author{\vspace{5pt} Eric A. Carlen$^{1}$,   Jos\'e A. Carrillo$^{2}$, and Michael Loss$^{3}$ \\
\vspace{5pt}\small{$1.$ Department of Mathematics, Hill Center, Rutgers University}\\[-6pt]
\small{
110 Frelinghuysen Road
Piscataway NJ 08854 USA}\\
\vspace{5pt}\small{$2.$ Instituci\'o Catalana de Recerca i Estudis Avan\c cats and Departament de Matem\`atiques}\\[-6pt]
\vspace{5pt}\small{Universitat Aut\`onoma de Barcelona, E-08193 Bellaterra, Spain}\\[-6pt]
\vspace{5pt}\small{$3.$ School of Mathematics, Georgia Institute of
Technology,} \\[-6pt]
\small{Atlanta, GA 30332 USA}\\
 }
\date{\version}
\maketitle \footnotetext [1]{Work of Eric Carlen and Michael Loss is partially supported by U.S.
National Science Foundation
grant DMS 0901632}
\maketitle \footnotetext [2]{Work of Jos\'e A. Carrillo is partially supported by the project MTM2008-06349-C03-03 DGI-MCI
(Spain) and 2009-SGR-345 from AGAUR-Generalitat de Catalunya.\\
\maketitle \footnotetext [3]Michael Loss  was partially supported by U.S.
National Science Foundation grant  DMS 0901304.

\copyright\, 2008 by the authors. This paper may be reproduced, in
its entirety, for non-commercial purposes.}

\begin{abstract}  We give a simple proof of the $\lambda = d-2$ cases of sharp Hardy-Littlewood-Sobolev inequalities for $d \geq 3$, and also the
sharp Log-HLS inequality for $d=2$,
via a monotone flow governed by the fast diffusion equation.

\end{abstract}

\noindent{\bf keywords:} Hardy-Littlewood-Sobolev | Fast Diffusion |
Gagliardo-Nirenberg-Sobolev.

\medskip

\section{Introduction}

We  explain an interesting relation between the
sharp Hardy-Littlewood-Sobolev (HLS) inequality for the resolvent of the Laplacian, the sharp
Gagliardo-Nirenberg-Sobolev (GNS) inequality, and the fast
diffusion equation (FDE). As a consequence of this relation, we
obtain a new  identity expressing the HLS functional as an integral
involving the fast diffusion flow and the GNS functional.  From this identity
we obtain a simple proof of the sharp HLS inequality in the cases
that express the regularizing properties of the Green's function
for the Laplacian in $\R^d$, for $d\ge3$, and of the Logarithmic
Hardy-Littlewood-Sobolev (Log-HLS) inequality, for $d=2$.
The proof also provides interesting information about the HLS functional  that
does not follow from previous proofs of the HLS inequality.

 Throughout the paper, we shall use
$\|f\|_p$ to denote the usual $L^p$ norms with respect to Lebesgue
measure:
$$
\|f\|_p = \left(\int_{\R^d}|f|^p\dd x\right)^{1/p}\, ,
$$
for $1 \leq p < \infty$.

\subsection{The sharp Hardy-Littlewood-Sobolev inequality}

The sharp form of the {\em HLS inequality} is due to Lieb
\cite{L83}. It states that for all non-negative measurable
functions $f$ on $\R^d$, and all $0 < \lambda < d$,
\begin{equation}\label{hls}
\frac{{\displaystyle\int_{\R^d}\int_{\R^d}
\frac{f(x)f(y)}{|x-y|^\lambda}\,\d x \,\d y}} {\|f\|_p^2 } \leq
\frac{{\displaystyle \int_{\R^d}\int_{\R^d}
\frac{h(x)h(y)}{|x-y|^\lambda}\,\d x \,\d y}} {\|h\|_p^2 }
\end{equation}
where
\begin{equation}\label{hdef}
h(x) = \left(\frac{1}{1+|x|^2}\right)^{(2d - \lambda)/2}\ .
\end{equation}
and
$p = 2d/(2d - \lambda)$.
Moreover, there is equality in(\ref{hls}) if and only if for some
$x_0\in \R^d$ and $s\in \R_+$, $f$ is a non-zero multiple of
$h(x/s-x_0)$.

The $\lambda = d-2$ cases of the sharp HLS inequality are
particularly interesting since they express the $L^p$ smoothing
properties of $(-\Delta)^{-1}$ on $\R^d$: for $d\geq 3$,
$$
\int_{\R^d} f(x)\left[ (-\Delta)^{-1} f\right](x)\dd x =
\frac{1}{(d-2)|S^{d-1}|}\int_{\R^d}\int_{\R^d}
\frac{f(x)f(y)}{|x-y|^{d-2}}\,\dd x\dd y\ ,
$$
where $|S^{d-1}|$ denotes the surface area of the $d-1$
dimensional unit sphere in $\R^d$. The integrals on the right hand
side of (\ref{hls}) can be computed explicitly in terms of
$\Gamma$-functions, and, after some computation with the
constants, one sees that for $\lambda =d-2$, (\ref{hls}) can be
rewritten as $\F[f]\geq 0$ for all $f\in L^{2d/(d+2)}(\R^d)$ where
\begin{equation}\label{hlsf}
{\mathcal F}[f] :=  C_S\|f\|^2_{2d/(d+2)} - \int_{\R^d} f(x)\left[
(-\Delta)^{-1} f\right](x)\dd x\ ,
\end{equation}
with
\begin{equation}\label{fdef2}
C_S :=  \frac{4}{d(d-2)}|S^d|^{-2/d}\ .
\end{equation}
We refer to this functional $\F$ on $L^{2d/(d+2)}(\R^d)$ as
the {\em HLS functional}.

Let $g$ be any smooth function of compact support, and let
$\langle g,f\rangle$ denote $\int_{\R^d}g(x)f(x)\dd x$. Then the
positivity of  $\F$ on $L^{2d/(d+2)}(\R^d)$ implies that for all
$f\in L^{2d/(d+2)}(\R^d)$,
$$
2\langle g,f\rangle - C_S\|f\|_{2d/(d+2)}^2 \leq  2\langle g,f\rangle -  \int_{\R^d} f(x)\left[
(-\Delta)^{-1} f\right](x)\dd x\ .
$$
Taking the supremum over $f$ on both sides; i.e., computing two
Legendre transforms, one finds
\begin{equation}\label{shso}
\frac{1}{C_S}\|g\|_{2d/(d-2)}^2 \leq \|\nabla g\|_2^2\ .
\end{equation}
Notice that $C_S$ is the least constant for which  (\ref{shso})
can hold for all smooth compactly supported functions $g$, since
the Legendre transform can be undone so that any improvement in
the constant in (\ref{shso}) would yield an improvement in the
constant in the HLS inequality, and this is impossible.

We can summarize the last paragraph by saying the that sharp HLS
inequality for $\lambda = d-2$, $d\geq 3$, is {\em dual} to the
sharp Sobolev inequality  (\ref{shso}), and because of this
duality, once one knows the sharp constant to one of these
inequalities, one knows the sharp constant to the other. A little
thought shows that the same is true for optimizers: Once one knows
the optimizers for one inequality, one knows the optimizers for
the other.

In this paper, we shall explain another kind of ``duality''
involving the $\lambda = d-2$ cases of the HLS inequality. This
duality relation, which does not have any evident connection with
the Legendre transform argument explained above,  relates the
$\lambda = d-2$ cases of the sharp HLS inequality to certain sharp
Gagliardo-Nirenberg-Sobolev  (GNS) inequalities, again with an
identification of sharp constants and optimizers. The GNS
inequalities in question are due, in their sharp form, to Del Pino
and Dolbeault \cite{DD}. They state that for all  locally
integrable functions $f$ on $\R^d$, $d\geq 2$, with a square
integrable distributional gradient, and $p$ with $1 < p <d/(d-2)$
\begin{equation}\label{gn}
\frac{ \|\nabla f\|_2^\theta \|f\|_{p+1}^{1-\theta}  }
{\|f\|_{2p}} \geq \frac{\|\nabla \tilde h\|_2^\theta
\|\tilde h\|_{p+1}^{1-\theta}}{\|\tilde h\|_{2p}}
\end{equation}
where
\begin{equation}\label{gdef}
\tilde h(x) = \left(\frac{1}{1+|x|^2}\right)^{1/(p-1)}\ .
\end{equation}
and
$\theta = d(p-1)/(p(d +2 - (d-2)p))$.
Moreover, there is equality in (\ref{gn}) if and only if for some
$x_0\in \R^d$ and $s\in \R_+$, $f$ is a non-zero multiple of
$\tilde h(x/s-x_0)$.

Notice that GNS optimizers are certain powers of  HLS optimizers,
and vice-versa. We shall see that this is no accident. In fact,
there is yet  another context in which the HLS optimizers play an
important role: they are the steady-state solutions of certain
nonlinear evolution equations pertaining to {\em fast diffusion}.

\subsection{The Fast Diffusion Equation}

The equation
\begin{equation}\label{fd0}
\frac{\partial u}{\partial t}(x,t) = \Delta u^m(x,t)
\end{equation}
with $0< m < 1$ describes fast diffusion. (For $m=1$, it is the
usual heat equation describing ordinary diffusion, and $m>1$
corresponds to {\em slow diffusion}.) As in \cite{cartos,CV}, note
that $u(x,t)$ solves (\ref{fd0}) if and only if
\begin{equation}\label{fd1}
v(x,t) := e^{td}u(e^t x,e^{\beta t})
\end{equation}
with $\beta=2- d(1-m)$ satisfies the
equation
\begin{equation}\label{fd}
\frac{\partial v}{\partial t}(x,t) = \beta \Delta v^m(x,t) +
\nabla \cdot[xv(x,t)]\ .
\end{equation}
For $m=1$, this is the Fokker-Planck equation, and (\ref{fd}) is a
non-linear relative of the Fokker-Planck equation. For all $1
-2/d< m < 1$, the equation (\ref{fd}) has integrable stationary
solutions. Computing them, one finds
\begin{equation}\label{fd2}
v_{\infty,M}(x) := \left(D(M) + \frac{1-m}{2\beta m}
|x|^2\right)^{-1/(1-m)}\ .
\end{equation}
The parameter $D(M)$ fixes the {\em mass } $M$ of the steady state
$v_{\infty,M}(x)$; i.e., the quantity
$$
M := \int_{\R^d} v_{\infty,M}(x) \, \d x \, .
$$
Computing the integral one finds that
$$
D(M) = C(d,m,\beta)M^{2/(2 - 2(1-m))}
$$
where $C(d,m,\beta)$ is a constant that may be expressed in terms
of $\Gamma$-functions.

The  self-similar solutions of (\ref{fd0}) corresponding through
the change of variables (\ref{fd1}) to the $v_{\infty,M}$ are
known as {\em Barenblatt solutions}, and $v_{\infty,M}$ is known
as the {\em Barenblatt profile} for (\ref{fd0}) with mass $M$.
Notice in the limiting case $m=1$, the Barenblatt profile
approaches a Gaussian, as one would expect. The Barenblatt
self-similar solutions are natural generalizations  of the
fundamental solutions of the heat equation. The Cauchy problem for
the FDE (\ref{fd0}) has been studied by many authors, we refer to
\cite{VaPME} for a full account of the literature.

It is established in  \cite{herpie} that the range of mass
conservation for the fast diffusion equation is $1-2/d < m < 1$.
As noted above, this is exactly the range of $m<1$ in which
integrable self-similar solutions exist. Within this range, the
flow associated to the fast diffusion equation is in many ways
{\em even better} than the flow associated to the heat equation;
see \cite{BV} and the references therein. The solutions of
(\ref{fd0}) with positive integrable initial data are $C^\infty$
and strictly positive everywhere instantaneously, just  as for the
heat flow.

Moreover, for non-negative initial data $f$ of mass $M$ satisfying
\begin{equation}\label{assump_inidat}
\sup_{|x|>R} f(x) |x|^{2/(1-m)} <\infty
\end{equation}
for some $R>0$, which means that $f$ is
decaying at infinity at least as fast as the
Barenblatt profile  $v_{\infty,M}$, the solution $v(x,t)$ of  (\ref{fd0})
with initial data $f$ satisfies the following remarkable bounds:
For any $t_*>0$, there exists a constant $C=C(t^*)>0$ such that
\begin{equation}\label{global bound}
\frac1C \leq \frac{v(x,t)}{v_{\infty,M}} \leq C \, ,
\end{equation}
for all $t\geq t_*$ and $x\in\R^d$. The lower bound is remarkable,
as our assumption on the initial data is an upper bound. This
shows ``how fast'' fast diffusion really is:  It spreads mass out
to infinity to produce the ``right tails'' instantly.

The proof of these bounds is based on the $L^\infty$-error
estimate obtained in \cite{Va} and improved to global Harnack
inequalities in \cite{BV}, see also \cite{CV}. Moreover, one can
show global smoothness estimates on the quotient, that is, for any
$t_*>0$
\begin{equation}\label{global bound reg}
\sup_{t\geq t_*}
\left\|\frac{v(\cdot,t)}{v_{\infty,M}}\right\|_{C^k(\R^d)} <\infty
\, ,
\end{equation}
for all $k\in\N$. Finally, it is well-known that
\begin{equation}\label{liml1}
\lim_{t\to\infty} \|v(t)-v_{\infty,M}\|_{L^1(\R^d)} =0 .
\end{equation}
For the best known rates of convergence see
\cite{BBDGV-arch}.

\section{Monotonicity of ${\mathcal F}$ along  fast diffusion}

Since the HLS minimizers are the attracting steady states for a
certain fast diffusion flow, one might hope that the HLS
functional ${\mathcal F}$ would be monotone decreasing along this
flow. This is indeed the case.

\begin{thm} \label{thm1}
Let $f \in L^{2d/(d+2)}(\R^d)$ be non-negative, and suppose that
$f$ satisfies (\ref{assump_inidat}) for some $R>0$, and $m=
d/(d+2)$, ensuring in particular that $f$ is integrable.  Let us
further suppose that
\begin{equation}\label{st1}
\int_{\R^d}f(x)\,\dd x =  \int_{\R^d}h(x)\,\dd x = M_*
\end{equation}
where $h$ is given by (\ref{hdef}) with $\lambda=d-2$. Let
$u(x,t)$ be the solution of (\ref{fd0}) with $m=d/(d+2)$ and
$u(x,1) = f(x)$. Then, for all $t>1$,
\begin{equation}\label{st2}
\frac{{\rm d}}{{\rm d}t}\F[u(\cdot,t)] = -2
\D[u^{(d-1)/(d+2)}(\cdot,t)] \leq 0
\end{equation}
where
\begin{equation}\label{ddef}
{\mathcal D}[g] \!:= \!
C_S\frac{d(d-2)}{(d-1)^2}\|g\|_{2d/(d-1)}^{4/(d-1)}\|\nabla
g\|_2^2 - \|g\|_{2(d+1)/(d-1)}^{2(d+1)/(d-1)}\ .
\end{equation}
\end{thm}

\noindent{\bf Proof:}  There are two things to be proved, namely the identity on the left hand side of  (\ref{st2}), and also the
non-negativity of the functional ${\mathcal D}$ defined in  (\ref{ddef}).

We begin with the former.  The uniform bounds on
the regularity of the quotient (\ref{global bound reg}) justify
all of the integration-by-parts used in the following computation of the derivative of $\F$
along the FDE flow for $m = d/(d+2)$:
\begin{equation}\label{dec4}
\frac{\partial}{\partial t}u(x,t) = \Delta u^{d/(d+2)}(x,t)\ .
\end{equation}
Therefore, let $u(x,t)$ solve (\ref{dec4}). We compute
\begin{align}\label{dec5}
\frac{\dd}{\dd t}\F[u] =&\, 2\,C_S\left(\int_{\R^d} u^{2m}\dd
x\right)^{2/d}\int_{\R^d}
 u^{(d-2)/(d+2)} \Delta u^m \,\dd x\nonumber\\
&- 2\int_{\R^d} \left(\Delta u^m\right) \left[(-\Delta)^{-1} u\right]\, \dd x\nonumber\\
=&\, -\frac{2C_S \,d(d-2)}{(d+2)^2}\left(\!\int_{\R^d}\!\!
u^{2m}\dd x\!\right)^{2/d} \!\!\!\!\int_{\R^d}\!\!
u^{-6/(d+2)}|\nabla u|^2\dd x\nonumber\\
&+ 2\int_{\R^d} u^{(2d+2)/(d+2)} u\dd x
\end{align}
Now define $g = u^{(d-1)/(d+2)}$. Then one computes
$$
\int_{\R^d}u^{-6/(d+2)}|\nabla u|^2\dd x = \left(\frac{d+2}{d-1}\right)^2\int_{\R^d}|\nabla g|^2\dd x\ .
$$

Rewriting the right hand side of (\ref{dec5}) in terms of $g$,
one finds
$$
\frac{\dd}{\dd t}\F[u] =\, -2\,C_S\left(\int_{\R^d}
g^{2d/(d-1)}\dd x\right)^{2/d}\frac{d(d-2)}{(d-1)^2}\int_{\R^d}
|\nabla g|^2 \dd x
+ 2\int_{\R^d} g^{(2d+2)/(d-1)}\dd x
\ .$$
Expressing this in terms of $\D[u]$, we have proved the left hand side of  (\ref{st2}).

We shall now show that $\D[u]$ is non-negative as a consequence of
the $p = (d+1)/(d-1)$ cases of the GNS inequality (\ref{gn}).
These can be written in the form
\begin{equation}\label{altgn}
C_{\rm GNS} \|\nabla g\|_2^2 \|g\|_{2d/(d-1)}^{4/(d-1)} \geq  \|g\|_{2(d+1)/(d-1)}^{2(d+1)/(d-1)}\ ,
\end{equation}
where, by definition, $C_{\rm GNS}$ is the best constant for which
this inequality is valid for all smooth $g$ with compact support.
One could compute the right hand side of  (\ref{gn}) to determine
the explicit value of  $C_{\rm GNS}$ and find that
\begin{equation}\label{conid}
C_{\rm GNS} = C_S\frac{d(d-2)}{(d-1)^2}\ .
\end{equation}

An easier way to see this is to go back to the first part of the
proof, and consider the initial data $f=h$, so that $u(x,t)$ does
not depend on $t$. Then by  what we just proved, ${\mathcal
D}(h^{(d-1)/(d+2)}) = 0$.  Notice that $h^{(d-1)/(d+2)}$ is an
optimizer for the $p = \frac{d+1}{d-1}$ case of  (\ref{gn}).
Hence, for the optimal $g$,
$$
C_S\frac{d(d-2)}{(d-1)^2} \|\nabla g\|_2^2 \|g\|_{2d/(d-1)}^{4/(d-1)} =  \|g\|_{2(d+1)/(d-1)}^{2(d+1)/(d-1)}\ ,
$$
and this proves  (\ref{conid}), and now the non-negativity of
${\mathcal D}$ follows from  (\ref{altgn}) and (\ref{conid}).
\lanbox

\

As we  show in the next section, all of the information that we
have used about the sharp GNS inequality can also be proved by a
fast diffusion flow argument without bringing anything else into
the argument.  Thus, while at present, our analysis may not
look self-contained, this will be remedied shortly. For now
though, let us finish with the HLS inequality.

As a direct consequence of the previous theorem, we deduce an
identity for the HLS functional that manifestly displays its
non-negativity.

\

\begin{thm} \label{thm2}
Let $f \in   L^{2d/(d+2)}(\R^d)$, $d\geq 3$ be non-negative.
Suppose also that $f$ satisfies
\begin{equation}\label{tails}
\sup_{|x|>R} f(x) |x|^{-(d+2)} <\infty
\end{equation}
for some $R>0$.
Then
\begin{equation}\label{st11}
\F[f]  =\frac{8}{d+2} \int_0^\infty\!\! e^{\beta t}
\D[u^{(d-1)/(d+2)}(\cdot,e^{\beta t})]\,\dd t \geq 0\ .
\end{equation}
Moreover, $\F[f]=0$ if and only if $f$ is a multiple of
$h(x/s-x_0)$ for some $s>0$ and $x_0\in \R^d$, with $h$ given by
(\ref{hdef}), $\lambda=d-2$.
\end{thm}
\medskip

\noindent{\bf Proof:} The assumption  (\ref{tails}) together with
the fact that $f\in L^{2d/(d+2)}(\R^d)$ implies the integrability
of $f$. Since  for all $\alpha>0$, $\F[\alpha f] = \alpha^2\F[f]$,
it is harmless to assume (\ref{st1}), which we do.  We may now
apply the previous theorem. Let $v(x,t)$ be the solution of
(\ref{fd}) with $v(x,0) = f(x)$. Let $u(x,t)$ be the solution of
(\ref{fd0}) with $u(x,1) = f(x)$. Because of the scaling relation
(\ref{fd1}), we have
\begin{equation*}
\F[v(\cdot,t)] = e^{t(d-2)}\F[u(\cdot, e^{\beta t})]\
\end{equation*}
with $\beta=4/(d+2)$. Then Theorem \ref{thm1} implies that, for
all $t>0$,
\begin{equation}\label{iden}
\frac{{\rm d}}{{\rm d}t}\left(e^{-t(d-2)} \F[v(\cdot,t)] \right) =
-\frac{8}{d+2} \, e^{\beta t} \D[u^{(d-1)/(d+2)}(\cdot,e^{\beta
t})] \ .
\end{equation}

We now claim that 
\begin{equation}\label{initval}
\lim_{t\to0}\F[v(\cdot,t)]= \F[f]\quad{\rm and}\quad 
\lim_{t\to\infty}\F[v(\cdot,t)]= 0\ .
\end{equation}
 Since $\F[h] = 0$, the latter fact  is an easy consequence of the 
the global bounds in (\ref{global bound}) due to
the assumptions (\ref{assump_inidat}) and (\ref{st1}), and a dominated convergence argument.  The former is slightly more subtle:
First, it is easy to show, using known facts about the Cauchy problem for the FDE \cite{VaPME}, that under our hypothesis,
$\lim_{t\to 0}\|v(\cdot,t) -f\|_{2d/(d+2)} = 0$. By an argument using Fatou's Lemma, the potential integral term can only jump downwards in the limit. 
Thus, at least we have $\F[f] \geq \lim_{t\to 0}\F[v(\cdot,t)]$,
and hence integrating  (\ref{iden}) over $[0,\infty)$, we obtain
$$
\F[f] \geq 
\frac{8}{d+2} \int_0^t e^{\beta s}
\D[u^{(d-1)/(d+2)}(\cdot,e^{\beta s})]\,\dd s \geq 0 \ .
$$
 
In particular, we have shown that $\F[f]\geq 0$
under the hypotheses of the theorem. Then an obvious truncation
and monotone convergence argument using the sequence of function $f_n = \min\{f,nh\}$
shows that $\F[f]$ is well defined, finite and non-negative for all non-negative 
$f\in L^{2d/(d+2)}(\R^d)$. 
This proves the $\lambda = d-2$ HLS inequality, and then by a standard argument using the positive definite nature of
the potential integral, shows that the potential integral is continuous on  $L^{2d/(d+2)}(\R^d)$. Thus,
$\F[f]$ is continuous on  $L^{2d/(d+2)}(\R^d)$, and   (\ref{initval}) now follows.
Now integrating  (\ref{iden}) over $[0,\infty)$ and using (\ref{initval}) yields the identity (\ref{st11}).

We now conclude from   (\ref{st11})  that $\F[f] = 0$ if and only if
$\D[u^{(d-1)/(d+2)}(\cdot,e^{\beta t})] = 0$ for all $t$. This is
equivalent to the existence of a constant $C$ and continuous
functions $s(t)$ and $x_0(t)$, defined for $t>1$ such that
$$
u(x,e^{\beta t})=Cs^{-d}(t)
\left[\tilde h\left(\frac{x}{s(t)}-x_0(t)\right)\right]^{(d+2)/(d-1)}
$$
due to the characterization of the optimizers in the GNS
inequality (\ref{gn}). Thus $u(x,e^{\beta t})$ is at each $t>0$ a
Barenblatt profile, thus by uniqueness of the Cauchy problem for
the FDE (\ref{fd0}), $u(x,e^{\beta t})$ is a self-similar
Barenblatt solution of the FDE (\ref{fd0}). Since $f(x) =
\lim_{t\to 0}u(x,e^{\beta t})$ in $L^1(\R^d)$, hence $f$ is itself
a Barenblatt profile, meaning that $f$ is a multiple of
$h(x/s-x_0)$ for some $s>0$ and some $x_0\in \R^d$. \lanbox

\

The identity   (\ref{st11}) has been derived under the hypothesis
(\ref{tails}). However, it is easy to  pass from
Theorem~\ref{thm2} to to the following, which is simply a
restatement of the $\lambda = d-2$ cases of Lieb's Theorem:

\begin{thm} \label{thm2a}
Let $f \in   L^{2d/(d+2)}(\R^d)$, $d\geq 3$ be non-negative. Then
$ {\mathcal F}[f]\geq 0$, and  $\F[f]=0$ if and only if $f$ is a
multiple of $h(x/s-x_0)$ for some $s>0$ and some $x_0\in \R^d$,
and where $h$ is given by (\ref{hdef}), $\lambda=d-2$.
\end{thm}

\

\noindent{\bf Proof:} We have already proved the inequality in the course of proving Theorem~\ref{thm2}.
 The cases of equality are somewhat more
subtle, and are dealt with in the next lemma.   \lanbox

\

\begin{lm}\label{lem1}
If $f\in L^{2d/(d+2)}(\R^d)$  is non-negative and satisfies
satisfies ${\mathcal F}[f]= 0$, then $f$ 
satisfies (\ref{assump_inidat}) for some $R>0$.
\end{lm}

To prove Lemma~\ref{lem1} we make use, for the first
time, of rearrangement inequalities  and the conformal invariance
of the HLS functional. 
Recently, Frank and Lieb gave  an interesting proof  of certain cases of the HLS inequality
\cite{FL} that uses  {\em reflection positivity} in place of rearrangements.  

\smallskip

\noindent{\bf Proof of Lemma~\ref{lem1}:}  
By a well--known theorem of Lieb
\cite{L77} on the cases of equality in the Riesz rearrangement
inequality, every optimizer $f$  in   $L^{2d/(d+2)}(\R^d)$ must be
a translate of its spherically symmetric decreasing rearrangement,
$f^*$. Making any necessary translation, we may assume $f=f^*$.
Next, as also shown by Lieb, the HLS functional is invariant under
the inversion mapping $f \mapsto \widehat f$ where
${\displaystyle
\widehat f(x) = |x|^{-(d+2)}f(x/|x|^2)}$,
which is an isometry on  $L^{2d/(d+2)}(\R^d)$.  Letting $x_0$ be
any unit vector,  $f$ is uniformly bounded on the unit ball
centered at $2x_0$. Thus $|x|^{-(d+2)}f(x/|x|^2-2x_0)$ is also an
optimizer, and satisfies (\ref{assump_inidat}) for some $R=1$. Now
the previous Theorem applies, and this function must be a
Barenblatt profile. It follows that the same is true of the
original $f$. \lanbox

\

Note that we have only used  the invariance of ${\mathcal F}$
under inversion, and hence under the full conformal group, to
settle the final points regarding cases of equality. It is
remarkable that  neither the fast diffusion flow, nor the GNS
inequalities possess this invariance, and yet  for a dense class
of functions  functions $f$, (\ref{st11}) expresses ${\mathcal
F}[f]$ in terms of the fast diffusion flow and  the GNS
functional.

\section{The sharp GNS inequalities and fast diffusion}

As we have seen, a calculation using fast diffusion reduces the proof of the 
$\lambda = d-2$ cases of the HLS inequality to the proof of certain GNS inequalities.
We now show, using results in \cite{cartos}, that another sort of 
calculation using a different  fast diffusion reduces the proof  these GNS inequalities to the
Schwarz inequality.

The FDE (\ref{fd}) with $1-2/d<m<1$ is a gradient flow of the
functional
\begin{equation}\label{entropy}
{\cal H}[v] = \int_{\R^d} \left[\frac{|x|^2}{2} v +
\frac{\beta}{m-1}v^m\right]\,dx \, ,
\end{equation}
with respect to the Euclidean Wasserstein distance, see
\cite{otto}. In particular, ${\cal H}[v]$ is a Liapunov functional
for (\ref{fd}), being its dissipation given by
\begin{equation}\label{disfast0}
\frac{d}{dt}{\cal H}[v] = - \int_{\R^d} \left| x+
\frac{m\beta}{m-1} \nabla v^{m-1}\right|^2 v\,dx :=-{\cal I}[v]
\end{equation}
for any solution $v(\cdot,t)$ to (\ref{fd}) and initial data
$v(x,0)$ satisfying the hypotheses of Theorem \ref{thm2}. Here,
the regularity properties of the solution (\ref{global bound}) and
(\ref{global bound reg}) that justified the computations in the
previous section ensure that at least when the initial data
satisfies  (\ref{assump_inidat}), the  dissipation  of the entropy
along the evolution is given by
\begin{equation}\label{disfast}
 \frac{d}{dt} {\cal I}[v] =  -2{\cal I}[v] -2(m-1)\int_{\R^d} v^{m} \left[ \Delta \xi \right]^2\,dx
 -2\int_{\R^d} v^m \left[ \sum_{i,j=1}^d
\left(\frac{\partial^2 \xi}{\partial x_i\partial x_j}\right)^2
\right] \, dx.
\end{equation}
with
$$
\xi=\frac{|x|^2}{2}+\frac{m\beta}{m-1} v^{m-1},
$$
as shown in  \cite{cartos}. Define
\begin{equation}\label{Rdef}
{\mathcal R}[v] :=  \int_{\R^d} v^{m} \left[ \Delta \xi \right]^2\,dx\ .
\end{equation}
By the Schwarz inequality for the Hilbert-Schmidt inner product,
$$
\int_{\R^d} v^m \left[ \sum_{i,j=1}^d \left(\frac{\partial^2
\xi}{\partial x_i\partial x_j}\right)^2 \right] \, dx\geq \frac1d
\int_{\R^d} v^{m} \left[ \Delta \xi \right]^2\,dx \, ,
$$
and thus from (\ref{disfast}) and  (\ref{Rdef}) we get
\begin{equation}\label{disfast2}
\frac{d}{dt} {\cal I}[v] \leq \,-2{\cal I}[v] -2 \left(m-1+\frac1d\right) {\cal R}[v]\ \, .
\end{equation}
As long as $(d-1)/d < m < 1$, the factor in from of ${\cal R}[v]$ is strictly positive.

Now  combine
(\ref{disfast0}) and (\ref{disfast2}) to conclude
\begin{equation}\label{disfast3}
\frac{d}{dt}{\cal H}[v] \geq \frac12 \frac{d}{dt} {\cal I}[v] +
\left(m-1+\frac1d\right){\cal R}[v]\, .
\end{equation}
Integrating this inequality in $t$  from $0$ to $\infty$,
and using the fact  that
$$
\lim_{t\to\infty}{\cal H}[v(\cdot,t)] = \lim_{t\to\infty}{\cal
I}[v(\cdot,t)] = 0,
$$
one gets
\begin{equation}\label{gnsct2}
{\cal H}[v(\cdot,0)] \leq  \frac12 {\cal I}[v(\cdot,0)] + \int_0^\infty \left(m-1+\frac1d\right){\cal R}[v(\cdot,t)]\dd t
\end{equation}
for all $v(\cdot,0)$ satisfying the assumptions of Theorem
\ref{thm2}. Since $\left(m-1+\frac1d\right){\cal R}[v] \geq 0$,
\begin{equation}\label{gnsct}
2{\cal H}[v(\cdot,0)] \leq {\cal I}[v(\cdot,0)]\ .
\end{equation}
This inequality is equivalent to the sharp GNS inequalities
(\ref{gn}). One see this by expanding the squares in this
inequality, cancelling the second moment terms from both sides,
and  performing an integration-by-parts allowed by (\ref{global
bound}). Then with  $m= (p+1)/(2p)$ and the change of dependent
variable $v(x,0)=:f(x)^{2p}$, and finally  using a standard
scaling argument, one arrives at (\ref{gn}); see \cite{DD} for
details. This finishes the summary of the relevant results in
\cite{cartos,DD}.

The exponent $m$ of the FDE
(\ref{fd0}) used to prove the particular
cases of the GNS inequalities involved in the proof of the HLS inequality in previous sections is $m=d/(d+1)$.

On the other hand,  the exponent $m$ in the FDE along which the HLS
functional is monotone is
$m=d/(d+2)$, which
corresponds to the critical exponent of FDE related to the
boundedness of the second moment of the stationary states
$v_{\infty,M}$, and it plays certain role in the large-time
assymptotics of the FDE, see \cite{CV,BBDGV-arch}.

We finally show how to extract from (\ref{disfast3}) the
characterization of the optimizers of the GNS inequalities
(\ref{gn}), at least under the conditions that are relevant for the application in the proof of Theorem~\ref{thm2}.

\begin{thm}\label{thm4}
Let $f$ be a positive measurable function on $\R^d$, $d\geq 2$,
with a square integrable distributional gradient $f$, such that
\begin{equation}\label{condopt}
\sup_{x\in\R^d} \frac{f^{-(p-1)/2p}(x)}{1+|x|^2} < \infty \, ,
\end{equation}
and $f$ being an optimizer of the GNS inequality (\ref{gn}). Then,
$f$ is given by (\ref{gdef}) up to translations and dilations.
\end{thm}

\

\noindent{\bf Proof:} Let us consider $v(x,0)=f^{2p}(x)$ as
initial data for the FDE (\ref{fd}) with $m = (p+1)/2p$.  Under
these conditions, we have derived (\ref{gnsct2}), and since
(\ref{gnsct}) is equivalent to the sharp GNS inequality for $f$
and $f$ is indeed a stationary state of the FDE (\ref{fd}) due to
(\ref{fd2}), it must be the case that
$$
\int_0^\infty {\cal R}[v(\cdot,t)]\,\d t = 0\ .
$$
Due to the positivity of $v(\cdot,t)$ for all $t>0$, we conclude
that $\Delta \xi = 0$ for all $t>0$.  It is straightforward to
infer from the global bounds (\ref{global bound}) that for any
$t_*>0$, there exists $D_1=D_1(t_*)>0$ such that
$$
\frac{1}{D_1} \leq \frac{|x|^2}{2}+\frac{m\beta}{m-1} v^{m-1}(x,t) \leq D_1
$$
for all $t\geq t_*$, $x\in\R^d$. Thus, $\xi$ is a globally bounded
harmonic function, and then Liouville's theorem implies that $\xi$
is constant. It follows that for each $t$, $v(\cdot,t)$ is a
Barenblatt profile, which determines the form of $f$ as in Theorem
\ref{thm2}. \lanbox

\section{Proof of the sharp Log-HLS inequality via fast diffusion}

In this section, we prove the sharp Log-HLS inequality on $\R^2$
by a similar fast diffusion flow argument. The sharp Log-HLS
inequality \cite{Be,CarLos} states that for all non-negative
measurable functions $f$ on $\R^2$ such that $f\ln f$ and $f\ln
(e+|x|^2)$ belong  to $L^1(\R^2)$,
\begin{equation}\label{loghls}
\int_{\R^2}\!\!\!f(x)\log f(x)\dd x+\!\frac{2}{M}\!
\int_{\R^2\times\R^2}\!\!\!\!\!\!\!\!\!\!\!\! f(x) \log|x-y|
f(y)\dd x \dd y\geq C\,,
\end{equation}
where $M:=\int_{\R^2}f\,\dd x$ with $C(M):=M(1+\log\pi-\log(M))$.
Moreover, there is equality if and only if $f(x) =
h_{\gamma,M}(x-x_0)$ for  some $\gamma > 0$ and some $x_0\in
\R^2$, where
\begin{equation}
  h_{\gamma,M}(x):=\frac{M}{\pi}\frac{\gamma}{{\left(\gamma +|x|^2
  \right)^2}}.
\end{equation}
Note that all of  integrals figuring in
the logarithmic HLS inequality are at least well defined with no
cancellation of infinities in their sum under the condition that
$f\ln f$ and $f\ln (e+|x|^2)$ belong to $L^1(\R^2)$.

We therefore define the Log HLS functional $\F$ by
\begin{equation}
\F[f] := \,\int_{\R^2}f(x)\log f(x)\dd x
+2\left( \int_{\R^2}f(x)\dd x\right)^{-1}\!\!\!\!
\iint_{\R^2\times\R^2}f(x) \log|x-y| f(y)\dd x \dd y \nonumber
\end{equation}
on the domain introduced above. The logarithmic HLS
functional involves three distinct integral functionals of $f$
while for $d\geq3$, the HLS functional involves only two. A more
significant difference is that the logarithmic HLS functional $\F$
is invariant under scale changes. That is, for $a>0$ and $f$ in
the domain of $\F$, define $f_{(a)} := a^2f(ax)$. One then
computes that $\F[f_{(a)}] = \F[f]$ for all $a>0$.

This simplifies our application of the fast diffusion equation, to
which we now turn. For $d=2$, $m = d/(d+2)$ reduces to $m = 1/2$.
Thus, the relevant cases of the fast diffusion equation in $d=2$
are
\begin{equation}\label{loghls3}
 \frac{\partial u}{\partial t}(x,t) = \Delta \sqrt{u}(x,t)\ ,
\end{equation}
and
\begin{equation}\label{loghls4}
 \frac{\partial v}{\partial t}(x,t) = \Delta \sqrt{v}(x,t) + \nabla \cdot[xv(x,t)]\ .
\end{equation}
As before, it is easily checked that the stationary states are
given by
$$
 v_{\infty,M}(x) = \left(D+\frac12 |x|^2\right)^{-2}\ ,
$$
for any mass $M>0$ with suitably chosen $D=D(M)$. Note that
\begin{equation}\label{loghls44}
h(x):=h_{1,4\pi}(x)=v_{\infty, M_*}(x)= \frac{4}{(1+|x|^2)^2}
\end{equation}
for a suitable $M_*$. For $d=2$, the scaling relation between
these two equations is that $u(x,t)$ solves (\ref{loghls3}) if and
only if $v(x,t) := e^{2t}u(e^tx,e^t)$ solves (\ref{loghls4}).
Notice that with $u$ and $v$ related in this way, the scale
invariance of $\F$ implies that
\begin{equation}\label{loghls55}
\F[v(\cdot,t)] = \F[u(\cdot,e^t)]\ .
\end{equation}

We now differentiate along the fast diffusion flow as before.  For
convenience, without loss of generality, we fix the initial mass.
We also impose the appropriate version of  (\ref{assump_inidat}).

\

\begin{thm}\label{thm3}
Let $f$ be a non-negative measurable functions  on $\R^2$ such
that $f\ln f$ and $f\ln (e+|x|^2)$ belong  to $L^1(\R^2)$. Suppose
that
$\int_{\R^2}f(x)\,\dd x =  \int_{\R^2}h(x)\,\dd x = M_*
$,
with $h$ given by (\ref{loghls44}). Them $\F[f] \geq \F[h]$,
and there is equality if and only if $f(x/s-x_0) = h(x)$ for some $s>0$ and some $x_0\in \R^2$.

Suppose in addition that
$f$ satisfies (\ref{assump_inidat}) for some $R>0$ and $m=1/2$.
Let $u(x,t)$ be the solution
of (\ref{loghls3}) with $u(x,1) = f(x)$. Then we have the identity
\begin{equation}\label{loghls8}
\F[f]  = \F[h] + \int_0^\infty \D[u^{1/4}(\cdot,e^t)]\,\dd t \geq
\F[h]\ ,
\end{equation}
where
$\D[g] \! = \|\nabla g\|_2^2\|g\|_4^4 - \pi \|g\|_6^6$
is non-negative by the $d=2$, $p=3$ case of the sharp GNS
inequality (\ref{gn}). 
\end{thm}

\

\noindent{\bf Proof:} Let 
 Let $v(x,t)$ be the solution
of (\ref{loghls4}) with $v(x,0) = f(x)$. Suppose initially that $f\leq Ch$ for some finite $C$. Under this additional argument
is it easy to prove  (\ref{initval}), though the $t=0$ limit is more subtle since in $d=2$, the potential integral is neither pointwise
positive, nor positive definite. However, our hypotheses ensure integrability of the positive and negative 
parts in the potential integral, and then monotone convergence may be used as before. A truncation argument, left to the reader,
then removes the extra assumption $f\leq Ch$.  

Differentiability of $\F[v(\cdot,t)] $ is justifies as before, and we have
$$
\F[h] - \F[f] = \int_0^\infty \frac{{\rm d}}{{\rm
d}t}\F[v(\cdot,t)] \,\dd t\ .
$$
But by (\ref{loghls55}),
${\displaystyle
\frac{{\rm d}}{{\rm d}t}\F[v(\cdot,t)] = \frac{{\rm d}}{{\rm d}t}\F[u(\cdot,e^t)]}$.
By the uniform regularity bounds on the
quotient (\ref{global bound reg}), we  compute 
\begin{align}
\frac{{\rm d}}{{\rm d}t}\F[u(\cdot,t)] = &\,-\frac{8\pi}{M_*}\int_{\R^2} \left[(-\Delta)^{-1}u\right] \Delta u^{1/2}\,\dd x \nonumber \\
 & \,+ \int_{\R^2}
 \log u \, \Delta u^{1/2}\d x \nonumber \\
= & \,\frac{8\pi}{M_*}\int_{\R^2}u^{3/2}\d x  -
\frac{1}{2}\int_{\R^2} \frac{|\nabla u|^2}{u^{3/2}}\d x\
.\label{redu}
\end{align}
Making the change of variables $g = u^{1/4}$, 
$\|u\|_{3/2}^{3/2} =\|g\|_6^6$, 
$\|g\|_4^4\d x = M_*$
and
$$
 \int_{\R^2}  \frac{|\nabla u|^2}{u^{3/2}}\d x = 16\|\nabla g\|_2^2\, ,
$$
leading together with (\ref{redu}) to (\ref{loghls8}). The proof of the statement about cases of equality proceed exactly as in 
Theorem \ref{thm2}.  The proof that  the condition
(\ref{assump_inidat}) may be relaxed as far as the inequality itself  is concerned, is similar,  with the
 the integrability condition on $f\ln (e+|x|^2)$ being used to insure integrability of the positive part of the potential integral. 
 \lanbox

\section{Consequences of the monotonicity}

The monotonicity of the HLS and Log-HLS functionals along fast
diffusion flows has interesting consequences. One of
these is the simplicity of the ``landscape''  of the 
Log-HLS  functional:
Let ${\mathcal C}$ be the (convex) set of non-negative functions on $\R^2$ satisfying all of the hypotheses
of Theorem~\ref{thm3}. 
Let ${\mathcal F}$ be  Log-HLS functional on
${\mathcal C}$. Then there are no strict local
minimizers of ${\mathcal F}$ in ${\mathcal C}$ other than the
absolute minimizers.

Indeed, this is an immediate consequence of
Theorems~\ref{thm1} and \ref{thm3}:  One can go monotonically down
to the absolute minimizers from any point in ${\mathcal C}$.
Of course, a similar result holds for the HLS funcional, but in this case the Euler-Lagrange equation is
thoroughly studied, and there are no non-negative critical points apart form global minimizers.

Next, as we have noted, the fast diffusion flow along which we
have shown ${\mathcal F}$, corresponding to the Log-HLS inequality
in $d=2$, to be monotone decreasing is gradient flow in the
$2$-Wasserstein metric for the entropy functional ${\mathcal H}$
defined in (\ref{entropy}). There is a kind of duality between the
HLS functional ${\mathcal F}$ and the entropy functional
${\mathcal H}$, as we now explain.

We first recall an observation of Matthes, McCann and Savare
\cite{MMS} concerning pairs of gradient flow equations. Consider
two smooth functions $\Phi$ and $\Psi$ on $\R^d$, and consider the
two ordinary differential equations describing gradient flow:
$$
\dot x(t) = -\nabla \Phi[x(t)]\qquad{\rm and}\qquad   \dot y(t) = -\nabla \Psi[y(t)]\ .
$$
Then of course $\Phi[(x(t)]$ and  $\Psi[(t(t)]$ are monotone
decreasing. Now differentiate each function along the others flow:
\begin{align*}
 \frac{{\rm d}}{{\rm d}t}\Phi[(y(t)] &= -\langle \nabla \Phi[y(t)],\nabla \Psi[y(t)]\rangle\\
  \frac{{\rm d}}{{\rm d}t}\Psi[(x(t)] &= -\langle \nabla \Psi[x(t)],\nabla \Phi[x(t)]\rangle\ .
\end{align*}
Thus, $\Phi$ is decreasing along the gradient flow of $\Psi$ for
any initial data if and only if $\Psi$ is decreasing along the
gradient flow of $\Phi$ for any initial data.

An analog of this holds for well-behaved gradient flows in the
$2$-Wasserstein sense, which is the result of \cite{MMS}. In our
case, we can apply it to the Log-HLS functional in $d=2$. Thus,
since ${\mathcal F}$ for the Log-HLS functional is decreasing
along the $2$-Wasserstein gradient flow for ${\mathcal H}$, one
can expect that ${\mathcal H}$ is decreasing along the
$2$-Wasserstein gradient flow for ${\mathcal F}$, which turns out
to be nothing other than the critical mass case of the
Patlak-Keller-Segel equation. Actually, the $m=1/2$, $d=2$ version
of ${\mathcal H}$ must be ``renormalized''  since in this case
$v_{\infty,M}$ does not have finite second moments, nor an
integrable square root. Still, this ``second Lyapunov function''
has been shown to be very useful in analyzing the critical mass
Patlak-Keller-Segel equation; see \cite{BCC}.

Finally, the main new results here, namely, the integral
identities (\ref{st11}) and (\ref{loghls8}), provide the starting
point of an analysis of ``remainder terms'' and ``stability''
results for the the sharp HLS and Log-HLS inequalities.  A quantitative stability theorem 
 shall be
developed elsewhere.  Finally, we expect to be able to carry out a similar proof for the cases $d-2 < \lambda < d$,
which involves a non-local analog of the fast diffusion equation.


\begin{thebibliography}{20}

\bibitem{Be}
W. Beckner, {\em Sharp {S}obolev inequalities on the sphere and
the {M}oser-{T}rudinger inequality}, Ann. of Math. 2, 138 (1993),
pp.~213--242.

\bibitem{BBDGV-arch}
A. Blanchet, M. Bonforte, J. Dolbeault, G. Grillo and J. L.
V\'azquez, {\em Asymptotic of the fast diffusion equation via
entropy estimates}, Arch. Rational. Mech. Anal., 191 (2009), pp.~
347--385.

\bibitem{BCC}
A. Blanchet, E. Carlen, and J.A. Carrillo, {\em Functional
inequalities, thick tails and asymptotics for the critical mass
Patlak-Keller-Segel model}, preprint.

\bibitem{BV}
M.~Bonforte and J.L.~V\'azquez, {\em Global positivity estimates
and Harnack inequalities for the fast diffusion equation}, Jour.
Func. Analysis, 240 (2006), pp.~399--428.

\bibitem{CS}
E.A.~Carlen and M.~ Loss, {\em Extremals of functionals with competing symmetries},
 J. Funct. Anal., 88 (1990), pp. 437--456

\bibitem{CarLos}
E.~Carlen and M.~Loss, {\em Competing symmetries, the logarithmic
{HLS} inequality and {O}nofri's inequality on {$S\sp n$}}, Geom.
Funct. Anal., 2 (1992), pp.~90--104.

\bibitem{cartos}
J.A.~Carrillo and G.~Toscani, {\em Asymptotic $L^1$ decay of the
porous medium equation to self-similarity}, Indiana Univ Math J.,
46 (2000), pp.~113-142.

\bibitem{CV}
J. A. Carrillo and J. L. V\'azquez, {\em Fine asymptotics for fast
diffusion equations}, Comm. Partial Differential Equations, 28
(2003), pp.~1023--1056.

\bibitem{DD}
M.~Del Pino and J.~Dolbeault, {\em Best constants for
Gagliardo-Nirenberg inequalities and applications to nonlinear
diffusions}, Jour. Math. Pures Appl., 81 (2002), pp.~301--342.

\bibitem{FL}
R.~Frank and E.H.~Lieb, {\em Inversion positivity and the sharp
HardyÐLittlewoodÐSobolev inequality}, to appear in Calculus of
Variations and PDEs.

\bibitem{herpie}
M.A.~Herrero and M.~Pierre, {\em The Cauchy problem for $u_t =
\Delta u^m$ when $0<m<1$}, Trans. AMS, 291 (1985), pp.~145--158.

\bibitem{L83}
E.H.~Lieb, {\em Sharp constants in the Hardy-Littlewood-Sobolev
and related inequalities}, Ann. of Math., 118 (1983),
pp.~349--374.

\bibitem{L77}
E.H.~Lieb, {\em Existence and uniqueness of the minimization solution of Choquard's non- linear equation},
Stud. Appl. Math., 57 (1977), pp.~93--105

\bibitem{MMS}
D. Matthes, R.~J. McCann, and G. Savar\'e, A family of nonlinear
fourth order equations of gradient flow type, Comm. Partial
Differential Equations, 34 (2009), pp.~1352--1397.

\bibitem{otto}
F. Otto, {\em The geometry of dissipative evolution equations: the
porous medium equation}, Comm. Partial Differential Equations, 26
(2001), pp.~101--174.

\bibitem{Va}
J. L. V\'azquez, {\em Asymptotic behaviour for the porous medium
equation posed in the whole space}, J. Evol. Equ., 3 (2003),
pp.~67--118.

\bibitem{VaPME}
J.L. V\'azquez, {\em The Porous Medium Equation. Mathematical
theory}, Oxford Mathematical Monographs, The Clarendon
Press/Oxford University Press, Oxford/New York, 2007.

\end{thebibliography}
\end{document}